\documentclass{article}
\usepackage{amssymb}
\usepackage{amsmath}
\usepackage{latexsym}
\usepackage{amsthm}
\usepackage{mathptmx}

\makeatletter
\def\thmhead@plain#1#2#3{%
  \thmname{#1}\thmnumber{\@ifnotempty{#1}{ }#2}%
  \thmnote{ \the\thm@notefont(#3)}}
\let\thmhead\thmhead@plain
\def\swappedhead#1#2#3{%
  \thmnumber{#2}\thmname{\@ifnotempty{#2}{. }#1}%
  \thmnote{ \the\thm@notefont(#3)}}
\makeatother

\newcounter{app}
\addtocounter{app}{1}

\theoremstyle{definition} 
\newtheorem{definition}{Definition}[app]

\theoremstyle{plain}      

\newcommand{\vold}[1]{\frac{d^3\vec{#1}}{2(2\pi)^{3/2} \omega(\vec{#1})}}

\begin{document}
\title{Deformation quantization of covariant fields}

\author{Giuseppe Dito\\
{\small {\it Laboratoire Gevrey de Math\'ematique Physique, CNRS \ UMR 5029}}\\
\vspace*{-1mm}
{\small{\it D\'epartement de Math\'ematiques, Universit\'e de Bourgogne}}\\
\vspace*{-1mm}
{\small{\it BP 47870, F-21078 Dijon Cedex, France.}}\\
\vspace*{-1mm}
{\small \texttt{Giuseppe.Dito@u-bourgogne.fr}}
}
\markboth{Giuseppe Dito}{Deformation quantization of covariant fields}
\maketitle

\begin{abstract}
After sketching recent advances and subtleties in classical
relativistically covariant field theories,
we give in this short Note some indications
as to how the deformation quantization approach can be used to
solve or at least give a better understanding of their quantization.
\end{abstract}


\section{Introduction}

We would like to discuss here a possible application of deformation 
quantization, viz. a nonperturbative construction of star-products (along 
the lines indicated in~\cite{Di93}), to the quantization of interacting 
field theory.  It is based on recent advances in the linearization programme 
for covariant wave equations proposed in~\cite{FS80a,FS80b} that uses 
nonlinear representation theory~\cite{FPS77}. We shall present this 
approach on an example. It can now be formulated in a more general setting 
thanks to results obtained by Flato, Simon, and Taflin concerning
the existence of global solutions, scattering theory,
and linearization of several covariant wave equations relevant 
to physics~\cite{ST85, FST87, ST93, FST97}
(see \cite{ST00} for a recent nontechnical survey of these results).

\section{A sketch of the linearization programme}

We shall give only the most basic facts and describe a few results,
relevant for our discussion of the linearization programme~\cite{FS80a,FS80b}
for covariant nonlinear equations. We skip the technical details: the reader 
can find them in the references listed at the end of the paper. At this stage
we shall only briefly present the context. For the reader's, convenience, 
we have collected in an appendix the basic notions on nonlinear
representation theory of Lie groups and algebras introduced and
developed in~\cite{FPS77} and some later works. 

Let $G$ be a Lie group and $E$ be some vector space of functions (endowed
with a Banach or Fr\'echet topology). We assume that there is an analytic 
action of $G$ on a neighborhood $\mathcal{O}$ of the origin in $E$ and that 
$0\in E$ is a fixed point for the action. Let $H(\mathcal{O},\mathcal{O})$
be the space of analytic maps from $\mathcal{O}$ to $\mathcal{O}$ 
{}\footnote{This kind of holomorphic mappings space between 
topological vector spaces can be given several topologies 
and the theory is  well developed, see~\cite{Din99}.}.
Hence there is a map  $g\in G \mapsto S_g \in H(\mathcal{O},\mathcal{O})$ such
that
$S_g(\mathcal{O})\subset \mathcal{O}$, $S_{g_1 g_2} = S_{g_1}S_{g_2}$,
$S_{g^{-1}}= S_g^{-1}$, and $S_g(0) =0$. The pair $(S,\mathcal{O})$
is called here an \textsf{analytic nonlinear representation} of $G$ in 
$\mathcal{O}$.

At the infinitesimal level, by differentiating the action of $S_{\exp(tX)}$,
where $X$ is an element of the Lie algebra $\mathfrak{g}$ of $G$,
one would expect to find ``vector fields'' $T_X$ on $\mathcal{O}$ 
(or on a open subset of it) satisfying $T_{[X,Y]} = T_X T_Y - T_Y T_X$.
The pair $(T,\mathcal{O})$ is called an \textsf{analytic nonlinear Lie algebra
representation} of $\mathfrak{g}$ on $\mathcal{O}$.

One may think of the maps $S$ and $T$ be given 
by their Taylor coefficients: $S=\sum_{n\geq1}S^n$ and $T=\sum_{n\geq1}T^n$.
Then the first order term $S^1$ (resp. $T^1$) is a linear representation
of $G$ (resp. $\mathfrak{g}$). Nonlinear representations
are built by successive extensions of a linear representation by its
successive tensor powers: they can be viewed as a kind of deformation of
their linear part. Interactions described by covariant nonlinear equations
(see below) can then be seen, as they should, as deformations of the 
free (noninteracting) situation, in line with Flato's deformation 
philosophy \cite{Fl82}.

\subsection{Covariant nonlinear equations}

We shall  be concerned with relativistic equations on Minkowski spacetime
with direct physical relevance and ``covariant'' will always mean covariant 
under the action of the Poincar\'e group 
$\mathcal{P}= \mathrm{SL}(2,\mathbb{C})\ltimes \mathbb{R}^4$.
The Lie algebra of $\mathcal{P}$ will be denoted by 
$\mathfrak{p} = \mathfrak{sl}(2,\mathbb{C})\dotplus \mathbb{R}^4$. 
We use standard conventions: Einstein summation
convention is used throughout for upper and lower indices 
(Latin indices run from $1$ to $3$ and Greek
ones from $0$ to $3$). The metric tensor $g_{\mu\nu}$ has signature $(+---)$. 
For convenience of notation we shall write $x\in\ \mathbb{R}^4$ 
as $x=(t,\vec{x})$, where $t=x^0$, $\vec{x} = (x^1,x^2,x^3)$. 
For a function $\phi\colon\mathbb{R}^4\rightarrow V$
where $V$ is a vector space, we set $\phi(t)\colon\mathbb{R}^3\rightarrow V$
by $\phi(t)(\vec{x})\equiv \phi(t,\vec{x}) \equiv \phi(x)$. An element 
$g\in \mathcal{P}$ is written as $g=(\Lambda,a)$ where 
$\Lambda\in \mathrm{SL}(2,\mathbb{C})$ and $a\in \mathbb{R}^4$. Let
$\Lambda \mapsto \tilde \Lambda$ be the canonical homomorphism 
of $\mathrm{SL}(2,\mathbb{C})$ onto the  group~$\mathrm{SO}_0(1,3)$.
Then the group law in $\mathcal{P}$ is 
$(\Lambda_1,a_1).(\Lambda_2,a_2) = (\Lambda_1\Lambda_2, a_1 
+ \tilde\Lambda_1 a_2)$.
The natural action
of the Poincar\'e group on spacetime will be denoted by 
$g\cdot x = \tilde \Lambda x + a$. Finally, the standard basis of 
$\mathfrak{p}$ is $\{P_\mu, M_{\alpha\beta}\}$ with $\alpha<\beta$.

Most of the classical wave equations met in physics can be cast into evolution
equations governed by some Hamiltonian operator~(see e.g. \cite{RS79}). 
Suppose that one starts with a vector space $E$ of initial data (usually a 
Banach or Fr\'echet space). For simplicity, elements of $E$ are functions 
defined on the whole of $\mathbb{R}^3$ and 
take their values in a finite-dimensional vector space $V$.
Let $H$ be an  operator on $E$ given by  a (non)linear 
differential operator such that $H(0)=0$. We are interested in the solutions
of the following evolution equation:
\begin{equation}\label{nle}
\frac{d\psi(t)}{dt} = H(\psi(t)),\qquad \psi(0)=\psi_0\in E.
\end{equation}
Now suppose that the previous evolution equation is associated to some
covariant nonlinear (wave) equation. So one has a finite-dimensional
representation of~$\mathrm{SL}(2,\mathbb{C})$ in~$V$, $\Lambda\mapsto A(\Lambda)$,
and given a solution $\psi$ of~Eq.~(\ref{nle}), then 
$\psi^g(x) \equiv A(\Lambda)\psi(\tilde\Lambda^{-1}( x -a))$ is
also a solution for any $g=(\Lambda,a)\in \mathcal{P}$. 
There is a local action $S$ of $\mathcal{P}$ on $E$ having the origin in $E$
as fixed point and given by 
$S_g(\psi_0)(\vec{x})=\psi^g(0,\vec{x})$, 
where $\psi_0(\vec{x}) = \psi(0,\vec{x})$, i.e., $S_{g_1g_2}=S_{g_1}S_{g_2}$.
In particular, if  $t\mapsto U_t$ is the evolution operator of~Eq.~(\ref{nle}), 
$U_{t+t'}=U_tU_{t'}$, $\psi(t)=U_t(\psi_0)$ and we have $U(t) = S_{\exp(tP_0)}$,
where $P_0\in \mathfrak{p}$ is the generator of time translations.
Hence we can associate to the covariant equation~(\ref{nle}) what looks very
much like a nonlinear representation of $\mathcal{P}$ in $E$. Conversely, if 
one has a nonlinear representation $(S,E)$ of $\mathcal{P}$, by formally 
differentiating $\psi(t)=S_{\exp(tP_0)}(\psi(0))$ with respect to $t$, one 
would expect to recover a covariant evolution equation as in~Eq.~(\ref{nle}).

The linearization programme~\cite{FS80a,FS80b} is based on the idea to take 
into account the full covariance of a nonlinear equation by associating to it 
a nonlinear representation, and then study its linearizability. In the case 
where the nonlinear  representation is linearizable, i.e., equivalent to its 
linear part which corresponds to the associated linear (free) equation, 
then one can directly deduce global existence in time of the solutions of the
nonlinear equation. 

As one would expect from deformation theory, the linearization involves the 
computation of some cohomology spaces. It is remarkable that 
formal nonlinear massive representations of the Poincar\'e group 
are always formally linearizable~\cite{Ta84}.
The covariance provides the right spaces where to
perform the study (spaces of differentiable or analytic vectors for the
linear part) and then to look for convergence of the expressions involved
when going to the analytic setting.

Notice that the Maxwell-Dirac system in $(1+3)$ dimensions 
(i.e., classical electrodynamics) presents many subtleties due to the 
zero mass of the photon. Even at the formal level, one cannot get a
(local) nonlinear representation of Lie algebra by using the techniques
of~\cite{FPS77}.  Further technical refinements are needed and it is a real 
\textit{tour de force} that was achieved in~\cite{FST97} by solving
the long standing problem of global solutions for classical electrodynamics and
asymptotic completeness. It turns out that Maxwell-Dirac system cannot be
linearized completely because of the long-range nature of the vector potential.
The asymptotic fields do not transform according
to a linear representation of the Poincar\'e group and consequently they
are not asymptotically free.
Nevertheless, one gets tractable expressions which 
provide new insights on the infrared problem in Quantum Electrodynamics.

\subsection{Nonlinear Klein-Gordon equation}\label{kg}

We  present an example that we shall need in 
Section~\ref{applidq} for our discussion on deformation quantization.
It is based on results obtained in~\cite{ST93}.
Consider the massive nonlinear Klein-Gordon equation:
\begin{equation}\label{nlkg}
\square \Phi  + m^2 \Phi = P(\Phi,\partial_t\Phi,\nabla\Phi), \quad m>0,
\end{equation}
where $\square=\partial^\mu\partial_\mu=\partial_t^2-\Delta$, 
and $P$ is a $C^\infty$ function covariant under $\mathcal{P}$ and 
vanishes at $0$ along with its first derivatives.

Simon and Taflin~\cite{ST85,ST93}  solved the Cauchy problem 
for~Eq.~(\ref{nlkg}),
in $(1+n)$-dimensional spacetime with $n\geq2$, by using linearization
techniques. They have shown global existence and asymptotic completeness
for Cauchy data in some neighborhood of the origin of an Hilbert space.

Equation~(\ref{nlkg}) can be put in the form  of an evolution 
equation~(\ref{nle}) by taking:
\begin{equation}\label{e1}
\psi=\begin{pmatrix}\Phi \\ \Pi \end{pmatrix}, \qquad
H(\psi) = \begin{pmatrix} 0 & I\\ \Delta -m^2 & 0 \end{pmatrix}
\begin{pmatrix}\Phi \\ \Pi \end{pmatrix} +
\begin{pmatrix}0\\ P(\Phi,\Pi,\nabla \Phi) \end{pmatrix}.
\end{equation}
We shall use the following decomposition into Fourier modes:
 \begin{equation}
  \begin{split}
 \Phi(t)(\vec{x})&=
   \int\vold{k} \big(\overline{a}(t)(\vec{k}) e^{-i\vec{k}\cdot \vec{x}}
      +a(t)(\vec{k})e^{i\vec{k}\cdot \vec{x}}\big), \\
 \Pi(t)(\vec{x})&=
  i \int \frac{d^3\vec{k}}{2(2\pi)^{3/2}}\big(\overline{a}(t)(\vec{k}) 
e^{-i\vec{k}\cdot \vec{x}} -a(t)(\vec{k})e^{i\vec{k}\cdot \vec{x}}\big),
  \end{split} \label{aat}
 \end{equation}
where $\omega(\vec{k})=\sqrt{|\vec{k}|^2 + m^2}$, along with:
 \begin{equation}\label{a+a-}
 a_+(t)(\vec{x})=
   i\int\frac{d^3\vec{k}}{(2\pi)^{3/2}} 
\overline{a}(t)(\vec{k}) e^{-i\vec{k}\cdot \vec{x}},\qquad
a_-(t)(\vec{x})=
-i \int \frac{d^3\vec{k}}{(2\pi)^{3/2}}a(t)(\vec{k})e^{i\vec{k}\cdot \vec{x}}.
\end{equation}
In terms of $\mathbf{a}=(a_+,a_-)$, the evolution equation for the 
nonlinear Klein-Gordon equation reads:
\begin{equation}\label{eva}
\frac{d}{dt}\begin{pmatrix}a_+(t) \\ a_-(t) \end{pmatrix} =
i\omega(-i\nabla)\begin{pmatrix}\phantom{-}a_+(t) \\ -a_-(t) \end{pmatrix}+
\begin{pmatrix}F(\mathbf{a}(t)) \\ F(\mathbf{a}(t)) \end{pmatrix},
\end{equation}
where $F(\mathbf{a}(t)) = P(\Phi(t),\Pi(t),\nabla \Phi(t))$ through the 
substitutions~(\ref{aat}) and (\ref{a+a-}) and\linebreak
 $\omega(-i\nabla)$ is the operator
acting as multiplication by $\omega(\vec{k})$ after Fourier transformation. 

The free Klein-Gordon equation (when $F=0$) induces a Lie algebra linear 
representation $T^1$ of $\mathfrak{p}$ in
$E_\infty=\mathcal{S}(\mathbb{R}^3,\mathbb{C})
\oplus \mathcal{S}(\mathbb{R}^3,\mathbb{C})$.
For 
$(f_+,f_-)\in E_\infty$ it is given by:
\begin{eqnarray*}
T^1_{P_0}(f_+,f_-)&=&i\omega(-i\nabla)(f_+, -f_-);\\
T^1_{P_j}(f_+,f_-)&=&\partial _j (f_+,f_-);\\
T^1_{M_{ij}}(f_+,f_-) &=& (x_i\partial _j - x_j \partial _i)(f_+, f_-);\\
T^1_{M_{0j}}(f_+,f_-)&=&i\omega(-i\nabla)( x_jf_+, - x_jf_-).
\end{eqnarray*}
The representation $T^1$ exponentiates to a continuous linear representation
of $\mathcal{P}$ in the Hilbert space 
$E=L^2(\mathbb{R}^3,\mathbb{C})\oplus L^2(\mathbb{R}^3,\mathbb{C})$,
for which $E_\infty$ is the set of $C^\infty$-vectors.

The nonlinear representation $T=\sum_{n\geq 1} T^n$ of $\mathfrak{p}$ 
in $E_\infty$
associated to Eq.~(\ref{eva}) is determined by the interaction term $F$. By
writing $T = T^1 +\tilde T$, on 
$\mathbf{f}=(f_+,f_-)\in E_\infty$, it is defined  by:
\begin{eqnarray*}
\tilde T_{P_0}(f_+,f_-)&=&(F(\mathbf{f}),F(\mathbf{f}));\\
\tilde T_{P_j}(f_+,f_-)&=&0;\\
\tilde T_{M_{ij}}(f_+,f_-) &=&0;\\
\tilde T_{M_{0j}}(f_+,f_-)&=&( x_jF(\mathbf{f}),x_jF(\mathbf{f})).
\end{eqnarray*}
The coefficients $T^n$, $n\geq2$, are the homogeneous terms of degree 
$n$ in $\tilde T$. With these notations, Eq.~(\ref{eva}) now reads:
$$
\frac{d\mathbf{f}(t)}{dt} = (T^1_{P_0} + \tilde T_{P_0})(\mathbf{f}(t)).
$$

It is shown in \cite{ST93} that this nonlinear representation of 
$\mathfrak{p}$ can be exponentiated to a nonlinear representation 
$U=\sum_{n\geq1}U^n$ of $\mathcal{P}$ 
on some neighborhood of $0\in E$ and that $U$ is linearizable. 
We shall be interested in the particular case where the interaction term~$P$
in~(\ref{nlkg}) is a covariant polynomial of the field $\Phi$ only, without 
constant term or terms of degree~$1$. Also we require that the classical
Hamiltonian associated to Eq.~(\ref{nlkg}) is a positive functional.
In that particular case, it follows directly from more general results
proved in~\cite{ST93} that:
\begin{itemize}
\item[(1)] $U$ is an analytic map 
$\mathcal{P}\times E_\infty \rightarrow E_\infty$
and $g\mapsto U_{g^{-1}}^1 U_g$ is a continuous map from $\mathcal{P}$
to the space of Banach analytic maps ${\cal H}(E_\infty,E_\infty)$.

\item[(2)] There exists a  unique invertible analytic 
$\Omega_\epsilon\colon E_\infty \rightarrow E_\infty$,
$\varepsilon=\pm$, such that $U_g\Omega_\varepsilon=\Omega_\varepsilon U_g^1$ 
for $g\in \mathcal{P}$.

\item[(3)] $\lim \limits_{t\to \pm\infty}
|| U_{\exp(tP_0)}(\Omega_\pm(\psi))
-U_{\exp(tP_0)}^1(\psi)||_E = 0, \quad \psi \in E_\infty$. 

\item[(4)] If the Cauchy data $(\Phi(0),\partial_t\Phi(0))$ are in
$E_\infty$, then $\Phi(t)$ is in $\mathcal{S}(\mathbb{R}^3,\mathbb{C})$
for all $t\in \mathbb{R}$ and $\Phi\in C^\infty(\mathbb{R}^4)$.

\end{itemize}

\section{Application of deformation quantization}\label{applidq}

Consider a finite-dimensional manifold $M$ with a Poison bracket 
$\{\cdot,\cdot\}$ defined on it. Recall that a deformation quantization 
of $M$ is given by a formal associative deformation 
$\star$ (star-product) of the usual multiplication of smooth functions on $M$:
\begin{equation}\label{3star}
f\star g = fg +\sum_{n\geq1} \lambda^n C_n(f,g),
\end{equation}
such that $[f,g]_\star = \frac{1}{2\lambda}(f\star g -g\star f)$ is a 
formal Lie algebra deformation of the Poisson bracket on $M$. 
All the series involved
are formal series in the deformation parameter $\lambda$ with coefficients
in the space of smooth functions on $M$ (see e.g.~\cite{DS} for details).
The cochains~$C_n$ are usually required to be bidifferential operators
vanishing on constants. In quantum-mechanical models, one takes 
$\lambda=\frac{i \hbar}{2}$
and the formal series often converge in the distribution sense.

In order to treat field theoretic problems, one would naturally
go over to infinite-dimensional manifolds of fields. Then one would consider
a family of functionals of the fields and try to deform it via a star-product.
Very quickly, one will be faced with the problem 
of giving a meaning to~(\ref{3star}) as a formal series, i.e., to have
well-defined cochains $C_n$ on the family of functionals considered. Of
course it is always possible to take a small enough family of functionals
such that the star-product is well defined. 
For example, if the manifold $M=H\oplus H$ where $H$ is
a Hilbert space, take as algebra of functionals $\mathcal{A}$
the one generated by (finite) tensor products of continuous linear forms on 
$H$, then it is easy to write down a well-defined star-product such that 
$\mathcal{A}\star\mathcal{A}\subset\mathcal{A}[\lambda]$. But this is of no
interest if one has in mind applications to field quantization: even
free Hamiltonians do not belong to this family of functionals. 
So our very first concern is to construct star-products defined on
a sufficiently rich family of functionals which should at least includes
the Poincar\'e generators.

For more than a half-century, physicists have known how to handle divergences
appearing in physical interacting quantum models, especially in QED,
and developed perturbative renormalization methods. The idea underlying most 
of the renormalization methods is to first regularize singular expressions 
(e.g. with a cut-off), then identify and extract terms becoming singular 
when the regularization is removed. This is done perturbatively with respect
to various relevant parameters of the model considered.

In the deformation quantization framework, one will obviously meet divergences
in the cochains of the star-product and try to remove them in a more 
controlled way. The idea here is to identify the diverging terms as singular 
cocycles or coboundaries of the Hochschild cohomology and get rid of them 
by passing to an equivalent star-product. This approach is perturbative 
in $\hbar$ and if one can manage to construct a ``renormalized''
star-product in such a way up to some order in $\hbar$, 
then it will necessarily be a multi-parameter deformation (the deformation 
parameters being the physical parameters of the model under study, e.g. 
coupling constant). The systematization of this idea would lead to a 
cohomological approach to renormalization (see~\cite{DS}, Section~4.2.1).

\subsection{Normal star-product}

Consider the free massive real scalar field $\Phi$. It satisfies
the linear Klein-Gordon equation $(\Box + m^2)\Phi(t,\vec{x})=0$, $m>0$.
Let $E_\infty=\mathcal{S}(\mathbb{R}^3,\mathbb{R})\oplus \mathcal{S}
(\mathbb{R}^3,\mathbb{R})$.
The initial data $(\phi,\pi)\in E_\infty$ are decomposed into 
(definite energy) Fourier modes by:
\begin{equation}\label{3fourier}
  \begin{split}
 \phi(\vec{x})&=
   \int\vold{k} \big(\overline{a}(\vec{k}) e^{-i\vec{k}\cdot \vec{x}}
      +a(\vec{k})e^{i\vec{k}\cdot \vec{x}}\big), \\
 \pi(\vec{x})&=
  i \int \frac{d^3\vec{k}}{2(2\pi)^{3/2}}\big(\overline{a}(\vec{k}) 
e^{-i\vec{k}\cdot \vec{x}} -a(\vec{k})e^{i\vec{k}\cdot \vec{x}}\big),
  \end{split} 
 \end{equation}
where as usual $\omega(\vec{k})=(|\vec{k}|^2 + m^2)^{1/2}$. 
Since $m>0$ and $(\phi,\pi)\in E_\infty$
it follows that $(\overline{a},a)$ is in 
$W_\infty\equiv\mathcal{S}(\mathbb{R}^3,\mathbb{C})
\oplus \mathcal{S}(\mathbb{R}^3,\mathbb{C})$.
We look at $W_\infty$ as a real vector space.

Denote by $\mathcal{A}=\mathcal{H}_s(W_\infty)$ the space of analytic 
$\mathbb{C}$-valued mappings on $W_\infty$ having semiregular kernels as 
tempered distributions. 
An element $F\in \mathcal{H}_s(W_\infty)$ can be written as:
\begin{equation}
F(\overline{a},a)
=\sum_{m,n\geq0} \langle F^{mn}, \overline{a}^{\otimes^m}
\otimes a^{\otimes^n}\rangle,
\end{equation}
where the kernel $F^{mn}\in \mathcal{S}'(\mathbb{R}^{3m},\mathbb{C})
\oplus \mathcal{S}'(\mathbb{R}^{3n},\mathbb{C})$ is regular with respect
to either its first $3m$ variables or its last $3n$ variables. Semiregularity
implies that usual functional derivatives of $F$ are in $\mathcal{S}$ as long 
as one does not consider mixed derivatives. Notice that the free Hamiltonian
$H_0(\overline{a},a)=\int \frac{d^3\vec{k}}{2} \overline{a}(\vec{k})a(\vec{k})$
belongs to $\mathcal{A}$ along with all of the free Poincar\'e generators.

We define normalized functional derivatives as:
\begin{equation*}
D_{a(\vec{k})} = (2\omega(\vec{k}))^{1/2} \frac{\delta}{\delta a(\vec{k})}
\qquad D_{{\overline{a}}(\vec{k})} = (2\omega(\vec{k}))^{1/2}
\frac{\delta}{\delta \overline{a}(\vec{k})},
\end{equation*}
and the standard Poisson bracket reads:
\begin{equation}\label{3poissonb}
\{F,G\} = 
\frac{2}{i}\int d^3 \vec{k} (D_{a(\vec{k})}(F)\ D_{{\overline{a}}
(\vec{k})}(G)\ - \ D_{{\overline{a}}(\vec{k})}(F)\  D_{a(\vec{k})}(G)\ ).
\end{equation}
Clearly $\mathcal{A}$ endowed with this bracket is a Poisson algebra 
(i.e., it is closed for the bracket operation), and
$(W_\infty, \mathcal{A}, \{\cdot,\cdot\})$ is a Poisson space.

The normal star-product $\star_N$ on $(W_\infty,\mathcal{A})$  is defined
by 
$$F\star_N G = F G + \sum_{n\geq1} \hbar^n C_n^N(F,G)$$
where:
\begin{equation}\label{norco}
C_n^N(F,G) = \frac{1}{n!} \int d^3 \vec{k_1}\cdots d^3 \vec{k_n}  
\big( D_{a(\vec{k_1})}\cdots  D_{a(\vec{k_n})}(F)\ 
D_{{\overline{a}}(\vec{k_1})}\cdots D_{{\overline{a}}(\vec{k_n})}(G)\big).
\end{equation}
Notice that the deformation parameter is $\frac{i \hbar}{2}$, the factor
$\frac{i}{2}$ has been absorbed in the definition of the cochains, so
the star-bracket is 
$[f,g]_{\star_N} = \frac{2}{i \hbar}(F\star_N G -G\star_N F)$.

The normal star-product realizes a deformation quantization of 
$(W_\infty, \mathcal{A}, \{\cdot,\cdot\})$ and it will be used in the next 
section as starting point for the construction of other star-products.

\subsection{Deformation quantization for interacting fields}

For our  discussion of the applications of the linearization programme
to deformation quantization of interacting fields, we shall consider
a massive real scalar field theory in $(1+3)$-dimensions with polynomial
interaction such that the Hamiltonian is positive:
\begin{equation}\label{lam4}
(\Box + m^2)\Phi(t,\vec{x})+V'(\Phi(t,\vec{x}))=0\qquad m>0,
\end{equation}
where $V'$ is the derivative of the potential in the Hamiltonian.
$V$ is assumed here to fulfill $V(0)=V'(0)=V''(0)=0$. We denote
by $\{P_\mu, M_{\alpha\beta}\}$ the Poincar\'e generators for
Eq.~(\ref{lam4}) seen as functionals of the Cauchy data $(\phi,\pi)$
or their Fourier modes $(\overline{a},a)$ defined above. 
Also $\{P^f_\mu, M^f_{\alpha\beta}\}$
will be the corresponding free generators.

For initial data $\Phi(0)=\phi$ and $\partial_t \Phi(0)= \pi$ in 
$E_\infty=\mathcal{S}(\mathbb{R}^3,\mathbb{R})\oplus \mathcal{S}
(\mathbb{R}^3,\mathbb{R})$,
it follows from the results of ~\cite{ST93}, reviewed in Section~\ref{kg},
that there exists a unique global solution $\Phi$
with free asymptotic fields $\Phi^\pm$:
$$
\lim_{t\rightarrow \pm\infty}\|\Phi(t) - \Phi^\pm(t)\|^{}_E =0,
$$
where $\| \cdot\|^{}_E$ denotes the Hilbertian norm of
$E=L^2(\mathbb{R}^3,\mathbb{R})\oplus L^2(\mathbb{R}^3,\mathbb{R})$.

Denote by $(\phi^\pm,\pi^\pm)$ the initial data of $\Phi^\pm$. 
The wave operators  $\Omega^\pm :
(\phi^\pm,\pi^\pm)\mapsto(\phi,\pi)$ are  invertible analytic 
mappings from  $E_\infty$ onto $E_\infty$,
they preserve the Poisson bracket (Poisson maps),
and  the classical scattering operator
$S = (\Omega^+)^{-1}\Omega^-$ is also analytic on $E_\infty$. 
Moreover, and more importantly for us, the wave operators linearize the 
Poincar\'e generators:
\begin{equation}\label{lin}
P_\mu\circ\Omega^{\pm}=P^f_\mu,\qquad 
M_{\alpha\beta}\circ\Omega^{\pm}=M^f_{\alpha\beta}.
\end{equation}

Notice that in general the generators $\{P_\mu, M_{\alpha\beta}\}$ do not
belong to the space $\mathcal{A}$ (e.g. for $\Phi^4$-theory). However
for the kind of potential we are considering, the first cochain
of the normal star-product $C^N_1$ makes still sense on the Poincar\'e 
generators, the remark applies to the Poisson bracket as well. 
However, the second cochain $C^N_2$
is not defined on the Poincar\'e generators $\{P_\mu, M_{\alpha\beta}\}$ 
in general.

We now introduce two star-products from the normal one:
\begin{equation}\label{starpm}
(F\star^\pm G)\circ\Omega^\pm = (F\circ\Omega^\pm)\star_N(G\circ\Omega^\pm).
\end{equation}
They are defined on $\mathcal{A}^\pm $ the push-forward of $\mathcal{A}$ by 
$\Omega^\pm$. Since $\Omega^\pm$ are Poisson maps the star-brackets 
of $\star^\pm$ are also deformations of the standard Poisson 
bracket~(\ref{3poissonb}). Since $\Omega^\pm$ depends (in a non trivial way!)
on the parameters (coupling constants) that may appear in the potential $V$ 
in~(\ref{lam4}), so will the cochains $C^\pm_n$ of $\star^\pm$. One should
also remark that (\ref{starpm}) does not necessarily imply the equivalence 
of the star-products.

{}From Eq.~(\ref{lin}), we deduce in particular for the $k$-th
$\star^\pm$-power of the interacting Hamiltonian $H=P_0$:
\begin{equation}\label{ham}
({\star^\pm } H)^k = ({\star^N } H_0)^k\circ(\Omega^\pm)^{-1},
\end{equation}
where $H_0= P_0^f$ is the free Hamiltonian. It is easy to check
that $({\star^N } H_0)^k$ belongs to $\mathcal{A}$ for all $k\geq1$,
and, since $(\Omega^\pm)^{-1}$ sends $W_\infty$ onto $W_\infty$, the right-hand
side of (\ref{ham}) is well defined for all $n\geq0$. Consequently all the 
powers $(\star^\pm H)^n$ are defined on $W_\infty$. Obviously, the same holds 
for any finite $\star^\pm$-products of the $\{P_\mu, M_{\alpha\beta}\}$. 
Therefore $\mathcal{A}^\pm $ contains a 
$\star^\pm$-subalgebra isomorphic to 
the universal enveloping algebra of the Poincar\'e Lie algebra.

A fundamental tool in deformation quantization is the star-exponential:
\begin{equation}\label{exp}
\mathrm{Exp}_\star(\frac{tF}{i\hbar})= 
\sum_{n\geq0} \frac{1}{n!}(\frac{t}{i\hbar})^n ({\star } F)^n.
\end{equation}
In many quantum-mechanical models, this series is convergent as a distribution,
and allows to determine the spectrum and eigenstates of some Hamiltonian or
any physically relevant observable. This is accomplished by performing
a Fourier analysis of $\mathrm{Exp}_\star(\frac{tH}{i\hbar})$ for some 
Hamiltonian. In field theory, for the free fields, one can easily sum up
the series and get expressions for relevant observables in a closed form.
For the example we are treating, each term in~(\ref{exp}) is well defined,
e.g. for the Hamiltonian.
However due to the presence of the factor $\frac{1}{\hbar}$ one cannot conclude
that the series is a formal series in both $\hbar$ and  $\hbar^{-1}$, because
their coefficients involve infinite sums. The linearization maps
would be useful for these considerations

The star-products $\star^\pm$ are best suited than the normal product as
they are free of a large class of divergences by construction. It is not clear
how the fields behave in the framework we have presented. A deeper
study of these constructions should open more avenues to treat divergences
problem in field theory and make more explicit the interplay between classical
and quantum field theories.

We end this Note by making some remarks about Electrodynamics. 
Maxwell-Dirac equations cannot be linearized and the asymptotic fields
are not free (the vector potential is free, but the Dirac spinor is not).
The existence of modified wave operators has been established in~\cite{FST97}
on a space of small initial data in neighborhood of the origin in some Hilbert 
space. These modified operators do linearize the Poincar\'e generators 
(seen as functional of the fields) and the electric current, but before
applying a procedure similar to the one we have sketched here, one will first
have to solve the deformation quantization of the asymptotic fields themselves.

\vskip5mm
\noindent\textbf{\large Acknowledgements}
\vskip2mm

\noindent The author thanks D. Sternheimer and E. Taflin for very
useful discussions and remarks.

\vskip5mm
\noindent\textbf{\large Appendix: Notions on nonlinear representations}
\vskip2mm

\noindent Let $E$ be a Fr\'echet space. Denote by $\mathcal{L}^n(E,E)$ the 
space of continuous symmetric $n$-linear maps on $E$ taking their values 
in $E$. The natural identification  
$\mathcal{L}^n(E,E)\sim \mathcal{L}(\hat{\otimes}^nE,E)$,
where $\hat\otimes$ is the completed symmetric projective tensor product, 
will be implicitly assumed. For $f\in \mathcal{L}^n(E,E)$,
we shall denote by $\hat f$ the associated homogeneous polynomial on $E$, i.e.,
$\hat f(x) = f(x,\ldots,x)$.  Let $\mathcal{F}(E)$ be the vector space
of formal series $F= \sum_{n\geq 1} f^n$, where $f^n\in \mathcal{L}^n(E,E)$.
There is an associative product on  $\mathcal{F}(E)$, denoted $\circ$, given
by the composition of formal series. For $F= \sum_{n\geq 1} f^n$ and 
$H= \sum_{n\geq 1} h^n$, the associated homogeneous polynomial of the
term of degree $n$ in $F\circ H$ reads:
\begin{equation}\label{circ}\tag{A.1}
(\widehat{F\circ H})^n (x) = \sum_{1\leq p\leq n} f^p\big(\sum_{n_1+\cdots 
+n_p=n}\hat{h}^{n_1}(x)\otimes\cdots\otimes\hat{h}^{n_p}(x)\big),\quad x\in E.
\end{equation}
Let $\mathcal{F}^\mathrm{inv}(E)$ be the group of invertible elements in 
$(\mathcal{F}(E),\circ)$.
\begin{definition}\label{SE}
A \textrm{formal nonlinear representation $(S,E)$ of a Lie group} $G$ in $E$ 
is a homomorphism $S\colon G \rightarrow \mathcal{F}^\mathrm{inv}(E)$ such that 
$(g;x_1,\ldots, x_n)\mapsto S_g^n(x_1,\ldots, x_n)$ are continuous  maps 
from $G\times E^n$ to $E$, $\forall n\geq1$.
\end{definition}
As a consequence of Eq.~(\ref{circ}), we have
$S^1_{gg'}(x) = S^1_g(S^1_{g'}(x))$ for $g,g'\in G$, $x\in E$. 
Thus the first order term of $S_g= \sum_{n\geq 1} S_g^n$ is a linear
representation of $G$ in $E$, denoted~$(S^1,E)$ and called the free part
of $(S,E)$.
\begin{definition}\label{TSE}
Two formal nonlinear representations $(S,E)$ and $(S',E)$ of a Lie\linebreak group~$G$
in $E$ are said to be formally equivalent, if there exists an element
$A\in \mathcal{F}^\mathrm{inv}(E)$ such that $S_g = A\circ S'_g\circ A^{-1}$, 
$\forall g\in G$. A formal nonlinear representation $(S,E)$ is called 
formally linearizable, if it is formally equivalent 
to the linear representation~$(S^1,E)$.
\end{definition}

Roughly speaking, a formal nonlinear representation can be seen as a formal 
action of~$G$ in~$E$ having $0\in E$ as a fixed point. By looking at the 
infinitesimal action, one would expect to encounter formal vector fields 
and their Lie bracket. It is indeed the case. For $F,H\in \mathcal{F}(E)$, 
define $F\bullet H \in \mathcal{F}(E)$ by:
\begin{equation}\label{bullet}\tag{A.2}
(\widehat{F\bullet H})^n(x)= 
\sum_{1\leq p\leq n} f^p\big(\sum_{1\leq i\leq p}
(x^{\otimes^{i-1}})\otimes \hat{h}^{n-p+1}(x)\otimes(x^{\otimes^{p-i}})\big),
\quad x\in E.
\end{equation}
The bracket $[F,H]\equiv F\bullet H - H \bullet F$ endows $\mathcal{F}(E)$
with a Lie algebra structure, and we have the corresponding of
Def.~\ref{SE} for
formal nonlinear representations of a Lie algebra.
\begin{definition}
A \textrm{formal nonlinear representation $(T,E)$ of a Lie algebra}
$\mathfrak{g}$ in $E$ is a Lie algebra homomorphism 
$T\colon \mathfrak{g} \rightarrow \mathcal{F}(E)$.
\end{definition}
By writing $T = \sum_{n\geq 1} T^n$, one can check that 
$T^1\colon \mathfrak{g}\rightarrow \mathcal{L}(E,E)$ is a linear Lie algebra
representation of $\mathfrak{g}$ in $E$, and $(T^1,E)$ is called the free
part of $(T,E)$.

There are \textrm{analytic} counterparts to the previous definitions. 
One would speak of analytic nonlinear representations of Lie groups and 
algebras, analytic equivalence of analytic nonlinear representations, etc. 
if the formal series involved are analytic in some neighborhood of the 
origin in $E$.

If one has a continuous linear representation $\pi$ of a real Lie group $G$ in 
a Banach space~$E$, then one has a linear representation $d\pi$ of the Lie 
algebra $\mathfrak{g}$ of $G$ in $E_\infty$, the Fr\'echet space of 
$C^\infty$-vectors of $\pi$. 
It is a remarkable fact that a similar statement holds  
for nonlinear representations. More precisely~\cite{FPS77, ST95}:
any (formal or analytic)  nonlinear representation $(S,E)$ of a real 
Lie group $G$ on a Banach space $E$ can be differentiated
and gives a nonlinear representation $(T,E_\infty)$ of the Lie algebra 
$\mathfrak{g}$  on the space of $C^\infty$-vectors for the free part 
$(S^1,E)$ of $(S,E)$. Then $T^1$ is equal to $dS^1$.
Conversely, if $G$ is moreover  connected and simply connected, any nonlinear
representation $(T,E_\infty)$ of $\mathfrak{g}$ such that its free part
$(T^1,E_\infty)$ is the differential of a linear representation $\pi$ of $G$
in a Banach space $E$, can be exponentiated to a unique nonlinear 
representation $(S,E_\infty)$ of $G$ such that $S^1=\pi$ (under some 
regularity assumptions on $(T,E_\infty)$, $(S,E_\infty)$ can be extended 
to a nonlinear representation $(S,E)$).


\begin{thebibliography}{FST97}

\bibitem[Din99]{Din99}
Dineen S.
\textit{Complex analysis on infinite-dimensional spaces}.
Springer-Verlag, London 1999.

\bibitem[Di90]{Di90}
Dito J. 
``Star-product approach to quantum field theory: the free scalar field,'' 
\textit{Lett. Math. Phys.} {\bf 20} (1990), 125--134;
 ``Star-products and nonstandard quantization for K-G equation,''
\textit{J. Math. Phys.} {\bf 33} (1992), 791-801; 
``An example of cancellation of infinities in star-quantization of fields,'' 
\textit{Lett. Math. Phys.} {\bf 27} (1993), 73--80.

\bibitem[Di93]{Di93}
Dito J.
\textit{Star-produits en dimension infinie~:
le cas de la th\'eorie quantique des champs},
Thesis, Universit\'e de Bourgogne, Dijon, 1993.

\bibitem[DS]{DS}
Dito G. and Sternheimer D.
``Deformation Quantization: 
Genesis, Developments and Metamorphoses ,'' In this volume.

\bibitem[Fl82]{Fl82}
Flato M. 
``Deformation view of physical theories," 
\textit{Czechoslovak J. Phys.} {\bf B32} (1982), 472--475. 

\bibitem[FPS77]{FPS77}
Flato M., Pinczon G. and Simon J.
``Non linear representations of Lie groups,''
\textit{Ann. Sci. \'Ecole Norm. Sup.} \textbf{10} (1977), 405--418.

\bibitem[FS80a]{FS80a}
Flato M. and Simon J.
``On a linearization program of nonlinear field equations,''
\textit{Phys. Lett. B} \textbf{94} (1980), 518--522.

\bibitem[FS80b]{FS80b}
Flato M. and Simon J.
``Linearization of relativistic nonlinear wave equations,''
\textit{J. Math. Phys.} \textbf{21} (1980), 913--917.

\bibitem[FST87]{FST87}
Flato M., Simon J. and Taflin E.
``On global solutions of the {M}axwell-{D}irac equations,''
\textit{Comm. Math. Phys.} \textbf{112} (1987) 21--49.


\bibitem[FST97]{FST97}
Flato M., Simon J.~C.~H. and Taflin E. 
{\it The Maxwell-Dirac equations: the Cauchy problem, asymptotic completeness 
and the infrared problem}, Mem. Amer. Math. Soc., 
{\bf 127} (number 606), 1997.


\bibitem[RS79]{RS79}
Reed M. and Simon B.
\textit{Methods of Modern Mathematical Physics III: Scattering Theory}.
Academic Press, Orlando 1979.

\bibitem[ST85]{ST85}
Simon J.C.H. and Taflin E.
``Wave operators and analytic solutions for systems of nonlinear
{K}lein-{G}ordon equations and of nonlinear {S}chr\"odinger equations,''
\textit{Comm. Math. Phys.} \textbf{99} (1985), 541--562.

\bibitem[ST93]{ST93}
Simon J.~C.~H. and Taflin E.
``The {C}auchy problem for nonlinear {K}lein-{G}ordon equations,''
\textit{Comm. Math. Phys.} \textbf{152} (1993), 433--478.

\bibitem[ST95]{ST95}  
Simon J.C.H.  and Taflin E.
``Initial data for non-linear evolution equations and differentiable 
vectors of group representations,''
In: J.~Bertrand et al.(eds.)
{\it Modern Group Theoretical Methods in Physics},
Math. Phys. Stud. {\bf 18},
Kluwer Acad. Publ., Dordrecht, 1995, 243--253.

\bibitem[ST00]{ST00}
Simon J.C.H. and Taflin E.
``Nonlinear relativistic evolution equations:
Survey of a new approach,''
In: G.~Dito and D.~Sternheimer (eds.)
\textit{Conf\'erence Mosh\'e Flato 1999}, 
Math. Phys. Stud. \textbf{21}, Kluwer Acad. Publ., Dordrecht,
2000, 389--402. 

\bibitem[Ta84]{Ta84}
Taflin, E.
``Formal linearization of nonlinear massive representations of the
connected Poincar\'e group,''
\textit{J. Math. Phys.} \textbf{25} (1984), 765--771.
\end{thebibliography}
\end{document}